\documentclass[graybox]{svmult}

\usepackage{mathptmx}       
\usepackage{helvet}         
\usepackage{courier}        
\usepackage{type1cm}        
\usepackage{makeidx}         
\usepackage{graphicx}        
\usepackage{multicol}        
\usepackage[bottom]{footmisc}

\usepackage{amssymb}
\usepackage{amsfonts}
 \newcommand{\nl}{\newline}
 \newcommand{\dist}{{\rm dist}}

\newcommand{\R}{{\mathbf{R}}}

 \newcommand{\diver}{{\rm div}}

 \newcommand{\half}{\frac{1}{2}}
 
 \newcommand{\darr}[4]{{\left\{\begin{array}{ll}
   {#1}&{#2}\\[0.2cm]
   {#3}&{#4}
 \end{array}\right.}}
  \newcommand{\darrsp}[4]{{\left\{\begin{array}{ll}
   {#1}&{#2}\\[0.4cm]
   {#3}&{#4}
 \end{array}\right.}}
 \newcommand{\tarr}[6]{{\left\{\begin{array}{lll}
   {#1}&{#2}\\
   {#3}&{#4}\\
   {#5}&{#6}
\end{array}\right.}}

\newcommand{\ia}{({\rm i})}
\newcommand{\ib}{({\rm ii})}
\newcommand{\ic}{({\rm iii})}

\newcommand{\parder}[2]{\frac{\partial{#1}}{\partial{#2}}}

\def\Frac{\displaystyle\frac}

\makeindex             

\begin{document}

\date{\em Dedicated to Ermanno Lanconelli on the occasion of his 70th birthday}

\title*{On the Hardy constant of some non-convex planar domains}
\titlerunning{On Hardy constant of planar domains} 
\author{ \bf Gerassimos Barbatis and Achilles Tertikas \\ [1.2 cm]
{\em Dedicated to Ermanno Lanconelli on the occasion of his 70th birthday}}
\authorrunning{G Barbatis \& A. Tertikas} 
\institute{Gerassimos Barbatis \at Department of Mathematics,
 University of Athens,  15784 Athens, Greece, \email{gbarbatis@math.uoa.gr}
\and Achilles Tertikas \at Department of Mathematics and Applied Mathematics,
 University of Crete, 70013 Heraklion, Greece and  \nl
Institute of  Applied and Computational Mathematics,
FORTH, 71110 Heraklion, Greece \email{tertikas@math.uoc.gr}}
%
%
\maketitle

\abstract*{The Hardy constant of a simply connected domain $\Omega\subset\R^2$ is the best constant for the inequality
\[
\int_{\Omega}|\nabla u|^2dx \geq c\int_{\Omega} \frac{u^2}{{\rm dist}(x,\partial\Omega)^2}\, dx \; , \;\;\quad  u\in C^{\infty}_c(\Omega).
\]
 After the work of Ancona where the universal lower bound 1/16 was obtained, there has been a substantial interest on computing or estimating the Hardy constant of planar domains. In \cite{BT} we have determined the Hardy constant of an arbitrary quadrilateral in the plane. In this work we continue our investigation and we compute the Hardy constant for other non-convex planar domains. In all cases the Hardy constant is related to that of a certain infinite sectorial region which has been studied by E.B. Davies.}

\abstract{The Hardy constant of a simply connected domain $\Omega\subset\R^2$ is the best constant for the inequality
\[
\int_{\Omega}|\nabla u|^2dx \geq c\int_{\Omega} \frac{u^2}{{\rm dist}(x,\partial\Omega)^2}\, dx \; , \;\;\quad  u\in C^{\infty}_c(\Omega).
\]
 After the work of Ancona where the universal lower bound 1/16 was obtained, there has been a substantial interest on computing or estimating the Hardy constant of planar domains. In \cite{BT} we have determined the Hardy constant of an arbitrary quadrilateral in the plane. In this work we continue our investigation and we compute the Hardy constant for other non-convex planar domains. In all cases the Hardy constant is related to that of a certain infinite sectorial region which has been studied by E.B. Davies.}

\vspace{11pt}

\vspace{6pt}

\section{Introduction}

The well-known Hardy inequality for $\R^N_+=\R^{N-1}\times (0,+\infty)$ reads
\begin{equation}
\int_{\R^N_+}|\nabla u|^2dx \geq \frac{1}{4}\int_{\R^N_+} \frac{u^2}{x_N^2}dx \; , \;\; \mbox{ for all }u\in C^{\infty}_c(\R^N_+),
\label{or1}
\end{equation}
where the constant $1/4$ is the best possible and equality is not attained in the appropriate Sobolev space.
The analogue of (\ref{or1}) for a domain $\Omega\subset\R^N$ is
\begin{equation}
\int_{\Omega}|\nabla u|^2dx \geq \frac{1}{4}\int_{\Omega} \frac{u^2}{d^2}\, dx \; , \;\; \mbox{ for all }u\in C^{\infty}_c(\Omega),
\label{or2}
\end{equation}
where $d=d(x)=\dist(x,\partial\Omega)$. However, (\ref{or2}) is not true without geometric assumptions on $\Omega$. The typical assumption made for the validity of (\ref{or2}) is that $\Omega$ is convex. A weaker geometric assumption introduced in \cite{bft} is that $\Omega$ is weakly mean convex, that is
\begin{equation}
-\Delta d(x)\geq 0 \; , \quad \mbox{ in }\Omega,
\label{c}
\end{equation}
where $\Delta d$ is to be understood in the distributional sense. Condition (\ref{c}) is equivalent to convexity when $N=2$ but strictly weaker than convexity when $N\geq 3$
\cite{ak}.
Other geometric assumptions on the domain that guarantee that the best Hardy constant is 1/4 were recently obtain in
\cite{av,Gk}.

For a general domain $\Omega$ we may still have a Hardy inequality provided that the boundary $\partial\Omega$ has some regularity. In particular it is well known that for any bounded Lipschitz domain $\Omega\subset\R^N$ there exists $c>0$ such that
\begin{equation}
\label{hohi}
\int_{\Omega}|\nabla u|^2 dx\geq c\int_{\Omega}\frac{u^2}{d^2}dx \;\; , \qquad \mbox{ for all }u\in C^{\infty}_c(\Omega).
\end{equation}
The best constant $c$ of inequality (\ref{hohi}) is called the Hardy constant of the domain $\Omega$.

In general the Hardy constant depends on the domain $\Omega$; see \cite{bl} for results that concern properties of this dependence. In dimension $N \geq 3$ Davies \cite{da1} has constructed Lipschitz domains with Hardy constant as small as one wishes. On the other hand for $N=2$ Ancona \cite{an} has proved that for a simply connected domain the Hardy constant is always at least $1/16$; see also \cite{laso} where further results in this directions where obtained.

Davies \cite{da1} computed the Hardy constant of an infinite planar sector $\Lambda_{\beta}$ of angle $\beta$,
\[ \Lambda_{\beta}  = \{ \; 0 < r , \;\;  0 < \theta  < \beta. \} \]
 He used the symmetry of the domain to reduce the computation to the study of a certain ODE; see (\ref{ode}) below. In particular he established the following two results, which are also valid for the circular sector of angle $\beta$:

(a) The Hardy constant is $1/4$ for all angles $\beta\leq \beta_{cr}$, where $\beta_{cr} \cong 1.546\pi$.

(b) For $\beta_{cr} \leq \beta\leq 2\pi$ the Hardy constant of $\Lambda_{\beta}$ strictly decreases with $\beta$ and at the limiting case $\beta=2\pi$ the Hardy constant is
$\cong 0.2054$.

Our interest is to determine the Hardy constant of certain domains in two space dimensions; see \cite{b,la} for relevant questions. In this direction, in our recent work \cite{BT} we have established

{\bf Theorem.} {\em 
Let $\Omega$ be a non-convex quadrilateral with non-convex angle $\pi<\beta <2\pi$. Then the Hardy constant of $\Omega$
depends only on $\beta$. The Hardy constant, which we denote from now on by $c_{\beta}$, is the unique solution of the equation
\begin{equation}
\label{tans1}
\sqrt{c_{\beta}}\tan\big( \sqrt{c_{\beta}}(\frac{\beta-\pi}{2})\big) = 2 \bigg( 
\frac{\Gamma(\frac{3+\sqrt{1-4c_{\beta}}}{4})}
{\Gamma(\frac{1+\sqrt{1-4c_{\beta}}}{4})}  
\bigg)^2,
\end{equation}
when $\beta_{cr}\leq\beta < 2\pi$ and $c_{\beta}=1/4$ when $\pi <\beta\leq \beta_{cr}$. The critical angle
$\beta_{cr}$ is the unique solution in $(\pi,2\pi)$ of the equation
\begin{equation}
\label{tans2}
\tan\big( \frac{\beta_{cr}-\pi}{4}\big) = 4\bigg(\frac{\Gamma(\frac{3}{4})}{\Gamma(\frac{1}{4})}\bigg)^2.
\end{equation}
}

Actually the constant $c_{\beta}$ coincides with the Hardy constant of the sector $\Lambda_{\beta}$, so equation (\ref{tans1}) provides an analytic description of the Hardy constant computed  numerically in \cite{da1}.

\begin{center}
\begin{figure}[here]
\begin{center}
\includegraphics[scale=0.26]{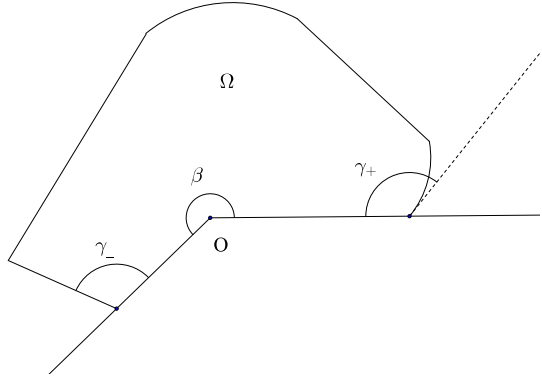}
\caption{A typical domain $\Omega$ for Theorem \ref{thm1}}
\label{fig:thm1}
\end{center}
\end{figure}
\end{center}

In this work we continue our investigation and determine the Hardy constant for other families of non-convex planar domains. Our first result reads as follows; see Fig.~\ref{fig:thm1}.

\begin{theorem}
Let $\Omega = K \cap\Lambda_{\beta}$, $\beta\in (\pi,2\pi]$, where $K$ is a bounded convex planar set and the vertex of $\Lambda_{\beta}$ is an interior point of $K$. Let $\gamma_+$ and $\gamma_-$ denote the interior angles of intersection of $K$ with $\Lambda_{\beta}$. There exists an angle $\gamma_{\beta}\in (\pi/2 ,\pi)$ such that if
$ \gamma_+ ,  \gamma_- \leq \gamma_{\beta}$, then the Hardy constant of $\Omega$ is $c_{\beta}$, where $c_{\beta}$
is given by (\ref{tans1}), (\ref{tans2}).
\label{thm1}
\end{theorem}
Detailed information on the angle $\gamma_{\beta}$ is given in Lemma \ref{lemma8} and Theorem \ref{new1}. We note that Theorem
\ref{thm1} can be extended to cover the case where $\Omega$ is unbounded and the boundary of the convex set $K$ does not intersect the boundary of the sector $\Lambda_{\beta}$; see Theorem \ref{unbdd}.

\begin{figure}[ht]
\begin{minipage}[b]{0.5\linewidth}
\centering
\includegraphics[scale=0.45]{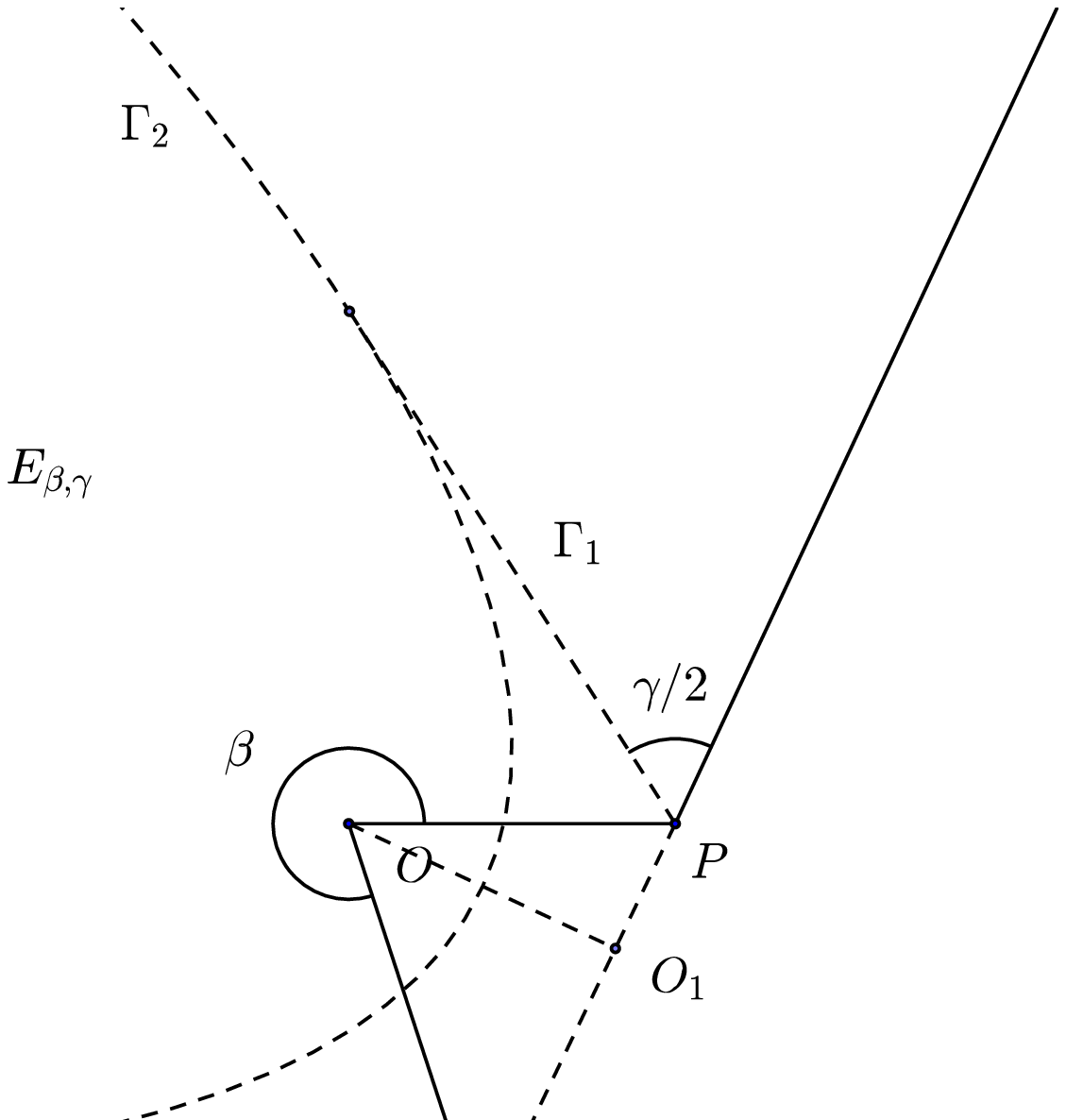}
\caption{A typical domain $E_{\beta,\gamma}$, $\gamma<\pi<\beta$}
\label{fig:ebg1}
\end{minipage}
\begin{minipage}[b]{0.5\linewidth}
\centering
\includegraphics[scale=0.2]{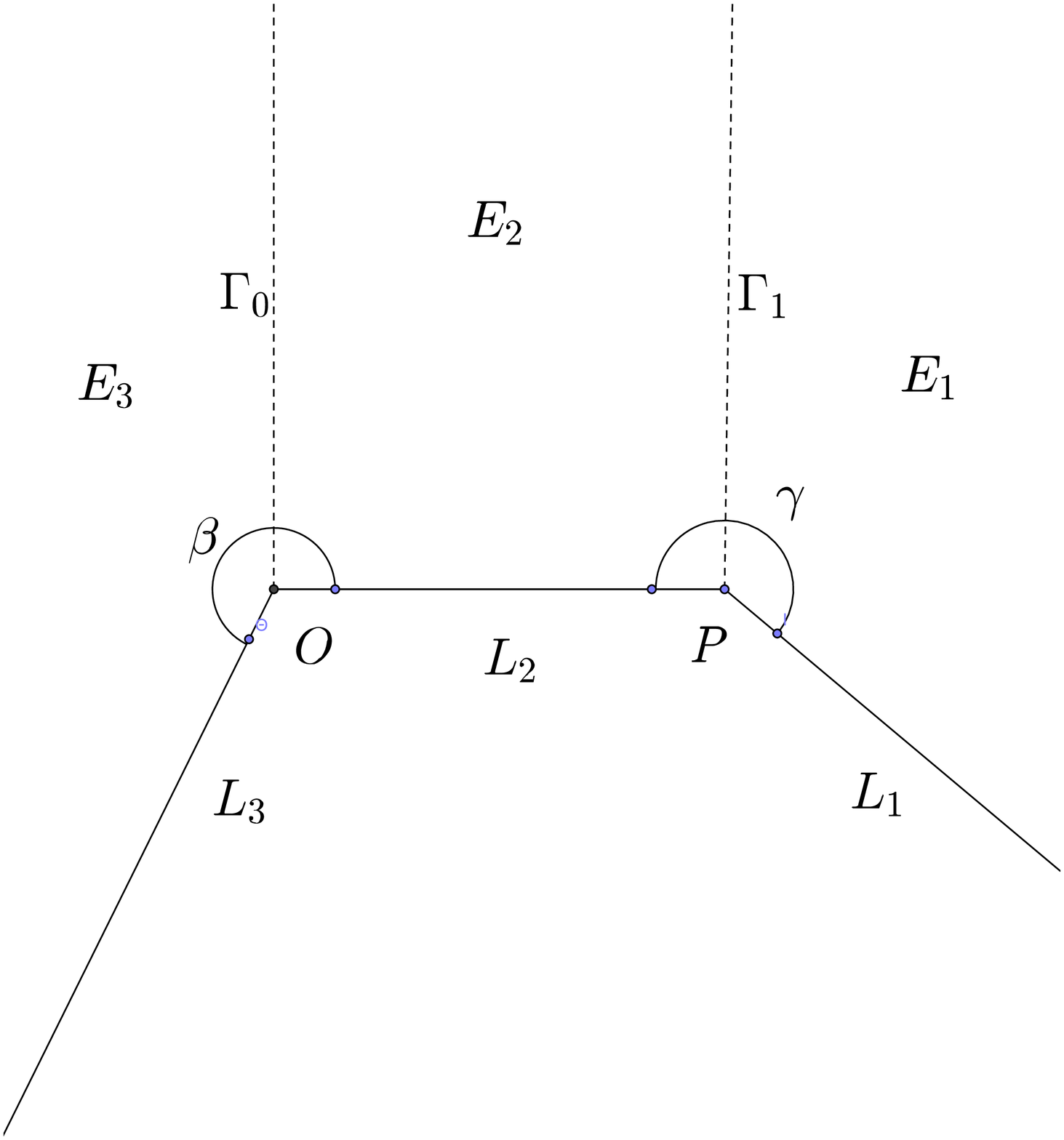}
\caption{A typical domain $E_{\beta,\gamma}$, $\beta , \gamma >\pi$}
\label{fig:ebg2}
\end{minipage}
\end{figure}

We next study the Hardy constant for a family of domains $E_{\beta,\gamma}$ which may have two non-convex angles.
The boundary $\partial E_{\beta,\gamma}$ of such a domain consists of the segment $ OP $ and two half lines starting from $ O $ and from $ P $ with  interior angles $\beta$ and $\gamma$; hence $\beta+\gamma\leq 3\pi$; see
Fig.~\ref{fig:ebg1} in case $\gamma<\pi$ and Fig.~\ref{fig:ebg2} in case $\gamma>\pi$. We then have

\begin{theorem}
\label{thm2}
$\ia$ If $0<\gamma \leq \pi \leq \beta \leq 2\pi$ then the Hardy constant of $E_{\beta,\gamma}$ is $c_{\beta}$.\nl
$\ib$ If $ \pi \leq \beta ,\gamma \leq 2\pi$ then the Hardy constant of $E_{\beta,\gamma}$ is $c_{\beta+\gamma -\pi}$, 
provided that
\begin{equation}
\label{oa}
|\beta-\gamma|  \leq \frac{2}{c_{\beta+\gamma-\pi}}\arccos(2\sqrt{c_{\beta+\gamma-\pi}}) .
\end{equation}
\end{theorem}
It is interesting to notice that in case (i) where we have only one non-convex angle, the Hardy constant is related to the non-convex angle $\beta$, whereas in case (ii) where we have two non-convex angles, the Hardy constant is related to the angle $\beta+\gamma-\pi$ formed by the two halflines.

\begin{center}
\begin{figure}[here]
\begin{center}
\includegraphics[scale=0.25]{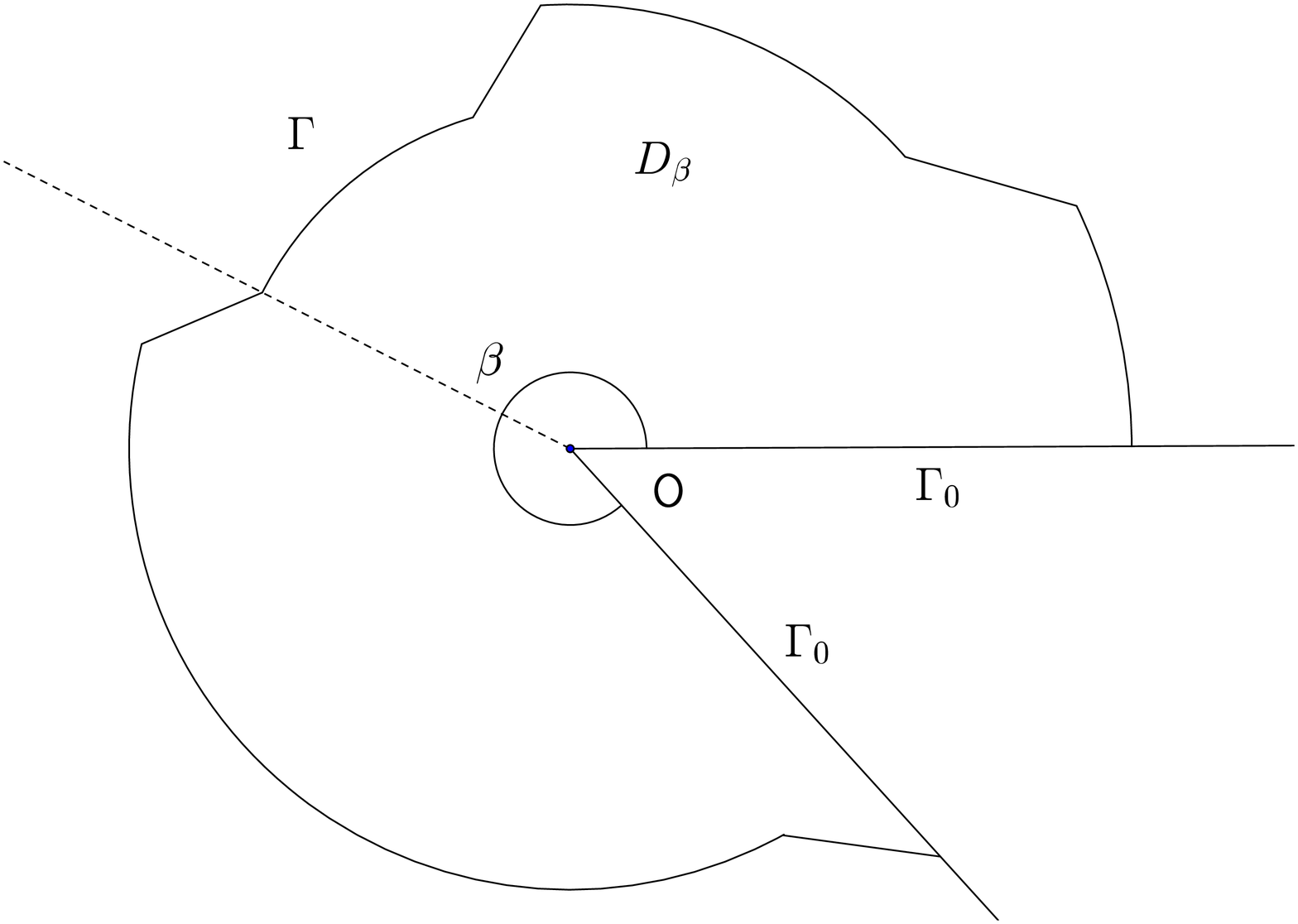}
\caption{A typical domain $D_{\beta}$. Note that $\Gamma$ is not necessarily the boundary of a convex set}
\label{fig:db}
\end{center}
\end{figure}
\end{center}

Our technique can actually be applied to establish best constant for Hardy inequality with mixed Dirichlet-Neumann boundary conditions.
We consider a bounded domain $D_{\beta}$ whose boundary $\partial D_{\beta}$ consists of two parts, $\partial D_{\beta}  =\Gamma_0 \cup \Gamma$. On $\Gamma_0$ we impose Dirichlet boundary conditions and it is from $\Gamma_0$ that we measure the distance from,
$d(x)={\rm dist}(x,\Gamma_0)$. On the remaining part $\Gamma$ we impose Neumann boundary conditions.
The curve $\Gamma_0$ is the union of two line segments which have as a common endpoint the origin $O$ where they meet at an angle $\beta$, $\pi<\beta\leq 2\pi$. 
We assume that the curve $\Gamma$ is the graph in polar coordinates of a Lipschitz function $r(\theta)$,
\[
\Gamma= \{ (r(\theta),\theta) :  0\leq \theta\leq \beta\} \, ;
\]
see Fig.~\ref{fig:db}.

We then have
\begin{theorem}
Let $D_{\beta}$ be as above, $\pi<\beta\leq 2\pi$. If $\Gamma$ is such that
\begin{eqnarray*}
&&  r'(\theta) \leq 0 , \quad 0\leq \theta\leq \frac{\beta}{2} , \\[0.2cm]
&&  r'(\theta) \geq 0 , \quad  \frac{\beta}{2}  \leq \theta\leq \beta \, ,
\end{eqnarray*}
then for all functions $u\in C^{\infty}(\overline{D_{\beta}})$ that vanish near $\Gamma_0$ there holds
\[
\int_{D_{\beta}}|\nabla u|^2 dx\, dy \geq c_{\beta}\int_{D_{\beta}}\frac{u^2}{d^2}dx\, dy \, .
\]
The constant $c_{\beta}$ is the best possible.
\label{thm3}
\end{theorem}

The structure of the paper is simple: in Section \ref{sec:lemmas} we prove various auxiliary results, while in Sections 3-5 we prove the theorems.

\section{Auxiliary results}
\label{sec:lemmas}

Let $\beta >\pi$ be fixed. We define the potential $V(\theta)$, $\theta\in (0,\beta)$,
\begin{equation}
V(\theta)=\tarr{\;\Frac{1}{\sin^2\theta} ,}{ 0< \theta  <\frac{\pi}{2} ,}{1,}{ \frac{\pi}{2}<\theta <\beta-\frac{\pi}{2},}
{\Frac{1}{\sin^2(\beta-\theta)},}{\beta-\frac{\pi}{2} <\theta <\beta.}
\label{v}
\end{equation}
For $c>0$ we then consider the following boundary-value problem:
\begin{equation}
\darr{ -\psi''(\theta) =c V(\theta)\psi(\theta),}{ 0\leq\theta\leq\beta ,}{ \psi(0)=\psi(\beta)=0\, }{}
\label{ode}
\end{equation}
It was proved in \cite{da1} that the Hardy constant of the sector $\Lambda_{\beta}$ coincides with
the largest positive constant $c$ for which (\ref{ode}) has a positive solution. Due to the symmetry of the potential $V(\theta)$ this also coincides with
the largest constant $c$ for which the following boundary value problem has a solution:
\begin{equation}
\darr{ -\psi''(\theta) =c V(\theta)\psi(\theta),}{ 0\leq\theta\leq\beta/2,}{ \psi(0)=\psi'(\beta/2)=0\, .}{}
\label{ode1}
\end{equation}
The largest angle $\beta_{cr}$ for which the Hardy constant is $1/4$ for
$\beta\in [\pi,\beta_{cr}]$ was computed numerically in \cite{da1} and analytically in \cite{BT,ti} where (\ref{tans2}) was established;
the approximate value is $\beta_{cr} \cong 1.546\pi$.

We define the hypergeometric function
\[
F(a,b,c ; z):=\frac{\Gamma(c)}{\Gamma(a)\Gamma(b)}\sum_{n=0}^{\infty}\frac{\Gamma(a+n)\Gamma(b+n)}{\Gamma(c+n)}\frac{z^n}{n!}.
\]
The boundary value problem (\ref{ode1}) was studied in \cite{BT} where the following lemma was proved:
\begin{lemma}
$\ia$ Let $\beta>\beta_{cr}$. The boundary value problem (\ref{ode1}) has a positive solution if and only if $c=c_{\beta}$. In this case the solution  is given by
\[
\psi(\theta)=\left\{
\begin{array}{l}
 \Frac{\sqrt{2}\cos\big(\sqrt{c}(\beta-\pi)/2\big)\sin^{\alpha}(\theta/2)\cos^{1-\alpha}(\theta/2) }{F(\half,\half, \alpha +\half ; \half)}F(\half,\half,\alpha +\half ; \sin^2(\frac{\theta}{2})), \\[0.2cm]
\hspace{8.1cm} \mbox{ if }\;0<\theta\leq \frac{\pi}{2}, \\[0.3cm]
\cos\big( \sqrt{c}(\frac{\beta}{2}-\theta) \big),  \hspace{5.45cm} \mbox{ if }\;\frac{\pi}{2} <\theta \leq \frac{\beta}{2}, 
\end{array}
\right.
\]
where $\alpha$ is the largest solution of $\alpha(1-\alpha)=c$.\nl
$\ib$ Let $\pi <\beta \leq \beta_{cr}$. The largest value of $c$ so that the boundary value problem (\ref{ode1}) has a positive solution is $c=1/4$. 
For $\beta=\beta_{cr}$ the solution is
\[
\psi(\theta)=\darrsp{ \Frac{\cos\big(\frac{\beta_{cr}-\pi}{4}\big)\sin^{1/2}\theta }{F(\half,\half, 1 ; \half)}F(\half,\half,  1; \sin^2(\frac{\theta}{2})),}
{0<\theta\leq \frac{\pi}{2},}{ \cos\big( \frac{1}{2}(\frac{\beta_{cr}}{2}-\theta) \big),}{ \frac{\pi}{2} <\theta \leq \frac{\beta_{cr}}{2}.}
\]
while for $\beta_{cr}<\beta<2\pi$ and $0<\theta<\pi/2$ it has the form
\begin{eqnarray*}
\psi(\theta) &=&c_1 \sin^{1/2}(\frac{\theta}{2}) \cos^{1/2}(\frac{\theta}{2})F(\frac{1}{2},\frac{1}{2},1;\sin^2(\frac{\theta}{2})) \\
&& + c_2 \sin^{1/2}(\frac{\theta}{2}) \cos^{1/2}(\frac{\theta}{2})F(\frac{1}{2},\frac{1}{2},1;\sin^2(\frac{\theta}{2}))
\int_{\sin^2(\theta/2)}^{1/2}\frac{dt}{t(1-t)F^2(\frac{1}{2},\frac{1}{2},1;t)}.
\end{eqnarray*}
for suitable $c_1$, $c_2$.
\label{lem:unify}
\end{lemma}

For our purposes it is useful to write the solution of (\ref{ode1}) in case $\beta\geq\beta_{cr}$
as a power series
\begin{equation}
\label{ps}
\psi(\theta) =\theta^{\alpha}\sum_{n=0}^{\infty}a_n\theta^{n} \; , 
\end{equation}
where $\alpha$ is the largest solution of the equation $\alpha(1-\alpha)=c_{\beta}$ in case $\beta>\beta_{cr}$ and $\alpha=1/2$ when $\beta=\beta_{cr}$.
We normalize the power series setting $a_0=1$; simple computations then give
\begin{equation}
a_1=0 \;\; , \qquad  a_2 =-\frac{\alpha(1-\alpha)}{6(1+2\alpha)}  .
\label{asymptotics}
\end{equation}
We also define the auxiliary functions
\begin{equation}
f(\theta) =\frac{\psi'(\theta)}{\psi(\theta)} \; , \qquad \theta\in (0, \beta) \; ,
\label{f}
\end{equation}
and
\begin{equation}
g(\theta) =\frac{\psi'(\theta)}{\psi(\theta)}\sin\theta \; , \qquad \theta\in (0, \beta) \; ,
\label{g}
\end{equation}
where $\psi$ is the normalized solution of (\ref{ode}) described in Lemma \ref{lem:unify}.
We note that these functions depend on $\beta$.
Simple computations show that they respectively solve the differential equations
\begin{equation}
f'(\theta) +f^2(\theta)+c_{\beta}V(\theta)=0 \; , \qquad 0<\theta <\beta
\label{def}
\end{equation}
and
\begin{equation}
\label{deg}
g'(\theta) =-\frac{1}{\sin\theta} \Big[ g(\theta)^2 -\cos\theta \, g(\theta) +c_{\beta} \Big] \; \; , \quad 0<\theta \leq\pi/2 .
\end{equation}

We shall also need the following
\begin{lemma}
Let $\pi\leq\beta\leq 2\pi$ and $\gamma\geq 0$ with $\beta+2\gamma\leq 3\pi$. Then
\[
f(\theta)\cos(\theta+\gamma)+\alpha[1+\sin(\theta+\gamma)]\geq 0 \; ,\quad 
\frac{\pi}{2}\leq\theta\leq\frac{3\pi}{2} -\gamma \, .
\]
\label{lem:quad10}
\end{lemma}
{\em Proof.} We first note that
\[
f(\theta)=\sqrt{c_{\beta}}\tan\Big( \sqrt{c_{\beta}}(\frac{\beta}{2}-\theta)\Big), \qquad\quad \frac{\pi}{2}\leq\theta\leq 
\frac{3\pi}{2}-\gamma,
\]
and
\[
-\frac{\pi}{4}\leq \sqrt{c_{\beta}}(\frac{\beta}{2}-\theta) \leq\frac{\pi}{4} \;, \qquad \quad \frac{\pi}{2}\leq\theta\leq \frac{3\pi}{2}-\gamma.
\]
It follows that the required inequality is written equivalently,
\begin{eqnarray}
\alpha(1&+&\sin(\gamma+\theta))\cos(\sqrt{c_{\beta}}( \frac{\beta}{2}-\theta)) \\ &+&\sqrt{c}\sin(\sqrt{c_{\beta}}( \frac{\beta}{2}-\theta))\cos(\gamma+\theta)\geq 0 \; , 
\;\; \frac{\pi}{2}\leq\theta\leq \frac{3\pi}{2}-\gamma.
\end{eqnarray}
But, since $\alpha\geq\sqrt{c_{\beta}}$,
\begin{eqnarray}
&\alpha&(1+\sin(\theta+\gamma))\cos(\sqrt{c_{\beta}}( \frac{\beta}{2}-\theta)) +\sqrt{c_{\beta}}\sin(\sqrt{c_{\beta}}( \frac{\beta}{2}-\theta))\cos(\theta+\gamma) \nonumber\\
&\geq&\sqrt{c_{\beta}}\Big\{(1+\sin(\theta+\gamma))\cos(\sqrt{c_{\beta}}( \frac{\beta}{2}-\theta)) +\sin(\sqrt{c_{\beta}}( \frac{\beta}{2}-\theta))\cos(\theta+\gamma)\Big\} \nonumber\\
&=& 2\sqrt{c_{\beta}} \sin\Big[ \sqrt{c_{\beta}}(\frac{\beta}{2}-\theta)+\frac{\pi}{4}+\frac{\theta}{2}+\frac{\gamma}{2}\Big] \sin(\frac{\pi}{4} +\frac{\theta}{2}+\frac{\gamma}{2}).
\label{s1}
\end{eqnarray}
The second sine is clearly non-negative, so it only remains to prove that the first sine is also non-negative. For this we use the monotonicity of 
$\sqrt{c_{\beta}}(\frac{\beta}{2}-\theta)+\frac{\pi}{4}+\frac{\theta}{2}+\frac{\gamma}{2}$ with respect to $\theta$ 
to obtain
\begin{eqnarray}
\sqrt{c_{\beta}}(\frac{\beta}{2}-\theta)+\frac{\pi}{4}+\frac{\theta}{2}+\frac{\gamma}{2}&\leq&
\sqrt{c_{\beta}}\big(\frac{\beta}{2}-(\frac{3\pi}{2}-\gamma)\big)+\frac{\pi}{4}+\frac{\frac{3\pi}{2}-\gamma}{2}+\frac{\gamma}{2} \nonumber \\ 
&=&\sqrt{c_{\beta}} \frac{\beta +2\gamma -3\pi}{2} +\pi \leq \pi ,
\end{eqnarray}
by our hypothesis $\beta +2\gamma\leq 3\pi$. This completes the proof. $\hfill\Box$

We shall need to consider the initial value problem (\ref{ivp}) below. Although this is a strongly singular problem, we shall see that standard comparison arguments hold. In particular we shall establish existence, uniqueness and monotonicity with respect to a parameter.
\begin{lemma}
Consider the singular initial value problem
\begin{equation}
\darr{ h'(\theta) =-\Frac{1}{\sin\theta} \Big( \alpha h(\theta)^2 -\cos\theta h(\theta) +1-\alpha \Big),}{
0<\theta\leq \frac{\pi}{2},}{h(0)=1.}{}
\label{ivp}
\end{equation}
$\ia$ If $\alpha \in (1/2, 1)$ then the problem has a classical solution which is unique. The solution $h(\alpha, \theta)$ depends monotonically on $\alpha$:
if $\alpha_1<\alpha_2$ then $h( \alpha_1, \theta) <  h( \alpha_2, \theta)$ for all $\theta\in (0,\pi/2]$.\nl
$\ib$ For $\alpha=1/2$ we do not have uniqueness. Indeed we have a continuum of positive solutions.\nl
$\ic$ Let $1/2<\alpha<1$ and in addition let $\overline{h}\in C[0,\pi/2] \cap C^1(0,\pi/2]$ be an upper solution of problem (\ref{ivp}), that is
\begin{equation}
\darr{ \overline{h}'(\theta) \geq -\Frac{1}{\sin\theta} \Big( \alpha \overline{h}(\theta)^2 -\cos\theta 
\overline{h}(\theta) +1-\alpha \Big) ,}{
0<\theta\leq \frac{\pi}{2},}{\overline{h}(0)\geq 1.}{}
\label{uivp}
\end{equation}
Then
\[
h(\alpha,\theta)\leq \overline{h}(\theta) \, , \qquad 0\leq\theta\leq\frac{\pi}{2} .
\]
\label{lem:h}
\end{lemma}
{\em Proof.} $\ia$ By Lemma \ref{lem:unify} the function
\[
\psi(\theta)=
\sin^{\alpha}(\theta/2)\cos^{1-\alpha}(\theta/2)F(\half,\half,\alpha +\half ; \sin^2(\frac{\theta}{2}))
\]
solves the differential equation
\[
\psi''(\theta) +\alpha(1-\alpha)\frac{\psi(\theta)}{\sin^2\theta} =0 \; , \quad  0<\theta<\frac{\pi}{2}.
\]
It is then easily verified that the function 
\[
h(\theta)=\frac{1}{\alpha}\frac{\psi'(\theta)}{\psi(\theta)}\sin\theta
\]
is a solution of the initial-value problem (\ref{ivp}).

We next establish the uniqueness of a solution. Let $h_1$, $h_2$ be two solutions of the initial value problem (\ref{ivp}). Then the function $z=h_2-h_1$ solves the singular linear initial value problem
\[
\darr{z'(\theta) = -\frac{1}{\sin\theta}\Big( \alpha(h_1+h_2) -\cos\theta\Big)z(\theta),}{}{z(0)=0.}{}
\]
Let us assume the $z$ is not identically zero. By the standard uniqueness theorem, $z$ cannot have any positive zeros, hence we may assume that $z(\theta)>0$ for all $\theta\in (0,\pi/2)$. However we have
$\alpha(h_1+h_2)-\cos\theta >0$ near $\theta=0$, hence $z$ decreases near zero, which is a contradiction.

The monotonicity of the solution $h$ with respect to $\alpha$ will follow from the monotonicity of the nonlinearity
with respect to $\alpha$.
Let
\[
V(\theta , h ,\alpha) =
-\frac{1}{\sin\theta} \Big( \alpha h^2 -\cos\theta h +1-\alpha \Big)
\]
For $0<h<1$ and $0<\theta<\pi/2$ we then have 
\begin{equation}
\frac{\partial V}{\partial\alpha} = \frac{1-h^2}{\sin\theta} >0.
\label{qqq}
\end{equation}
Now, let $1/2<\alpha_1<\alpha_2<1$.
By (\ref{qqq}) we have $h(\alpha_2,\theta)>h(\alpha_1,\theta)$ near $\theta=0$. Once we are away from $\theta=0$ we can apply the standard comparison arguments to complete the proof.

$\ib$ By Lemma \ref{lem:unify} 
the general solution of the equation
\[
\psi''(\theta) +\frac{1}{4}\frac{\psi(\theta)}{\sin^2\theta} =0 \; , \quad  0<\theta<\frac{\pi}{2} ,
\]
is
\begin{eqnarray*}
\psi(\theta) &=&c_1 \sin^{1/2}(\frac{\theta}{2}) \cos^{1/2}(\frac{\theta}{2})F(\frac{1}{2},\frac{1}{2},1;\sin^2(\frac{\theta}{2})) \\
&+& c_2 \sin^{1/2}(\frac{\theta}{2}) \cos^{1/2}(\frac{\theta}{2})F(\frac{1}{2},\frac{1}{2},1;\sin^2(\frac{\theta}{2}))
\int_{\sin^2(\theta/2)}^{1/2}\frac{dt}{t(1-t)F^2(\frac{1}{2},\frac{1}{2},1;t)}.
\end{eqnarray*}
This is positive in $(0,\pi/2]$ when $c_1>0$ and $c_2\geq 0$. 
For any such $\psi$ the function
\[
h(\theta)=\frac{2\psi'(\theta)}{\psi(\theta)}\sin\theta
\]
then satisfies
\[
h'(\theta) = -\frac{1}{2\sin\theta} \Big( h(\theta)^2 -2\cos\theta h(\theta) +1\Big) \; , \;\;\; h(0)=1.
\]
Actually after some computations we find that the function $h$ is given in this case by
\begin{eqnarray*}
h(\theta)&=&\cos\theta +\sin^2\theta \frac{F(\frac{3}{2},\frac{3}{2},2;\sin^2(\frac{\theta}{2}))}{4F(\frac{1}{2},\frac{1}{2},1;\sin^2(\frac{\theta}{2}))} \\
&-& \frac{\lambda}{F^2(\frac{1}{2},\frac{1}{2},1;\sin^2(\frac{\theta}{2}))
\Big(  
1+\lambda\int_{\sin^2(\theta/2)}^{1/2}\frac{dt}{t(1-t)F^2(\frac{1}{2},\frac{1}{2},1;t)}
\Big)},
\end{eqnarray*}
where $\lambda =c_2/c_1\geq 0$.

$\ic$ When $\overline{h}(0)>1$ the result follows immediately by combining continuity with standard comparison arguments. Assume now that $h(0)=1$. The function $z=\overline{h}-h$ then satisfies
\begin{equation}
\darr{z'(\theta) \geq -\frac{1}{\sin\theta}\Big( \alpha(\overline{h} +h) -\cos\theta\Big)z(\theta),}{}{z(0)=0.}{}
\label{lll}
\end{equation}
The quantity $\alpha(\overline{h} +h) -\cos\theta$ is positive near $\theta=0$, say in $(0,\theta_0)$. We shall establish that $z\geq 0$ in this interval; the result for $(0,\pi/2)$ will then follow immediately.
Suppose on the contrary that there exists an interval $(\theta_1,\theta_2)\subset (0,\theta_0)$ such that $z<0$ in
$(\theta_1,\theta_2)$. By (\ref{lll}) we conclude that $z$ is actually strictly increasing in $(\theta_1,\theta_2)$. This contradicts the initial value $z(0)=0$.
$\hfill\Box$

From Lemma \ref{lem:h} it follows that the case $\alpha=1/2$ is critical and needs a different approach. This will be done in the next lemma. 
In order to make explicit the dependence on $\beta$ we denote
\[
g(\beta,\theta)=\frac{\psi_{\theta}(\beta,\theta)}{\psi(\beta,\theta)}\sin\theta \; ,
\]
where $ \psi(\beta,\theta) $ is the solution of (\ref{ode}) and $ \psi_{\theta}(\beta,\theta) $ is the derivative with respect to $ \theta $.
\begin{lemma}
\label{222aa}
Suppose $\pi\leq\beta\leq \beta_{cr}$. Then $g(\beta,\theta)$, $0 < \theta\leq\pi/2$, is strictly increasing as a function of $\beta$, that is, if $\pi\leq\beta_1<\beta_2\leq \beta_{cr}$ then
$g(\beta_1,\theta) <  g(\beta_2,\theta)$ for all $\theta\in (0,\pi/2]$.
\end{lemma}
{\em Proof.}
The function $g(\beta,\theta)$ solves the differential equation
\begin{equation}
\label{eqg}
\parder{g}{\theta} = -\frac{1}{\sin\theta}\Big( g^2 -g \cos\theta  +\frac{1}{4} \Big).
\end{equation}
Since
\[
g(\beta,\frac{\pi}{2})= \frac{1}{2}\tan(\frac{\beta-\pi}{4}),
\]
which is strictly increasing with respect to $\beta$, the result follows from a standard comparison argument.
$\hfill\Box$

Let us note here that for $\pi\leq\beta\leq\beta_{cr}$ we have $g(\beta,0)=1/2$. So the functions $g(\beta,\cdot)$,
$\pi\leq\beta\leq\beta_{cr}$, all solve the same initial value problem.

\begin{lemma}
Let $\beta\in [\pi,2\pi]$. There exists an angle $\gamma^*_{\beta}$ so that for all $0 <\gamma\leq\gamma^*_{\beta}$ we have
\begin{equation}
g(\beta,\theta)\cos(\theta +\frac{\gamma}{2}) + \alpha\cos\frac{\gamma}{2}\geq 0 \; , \qquad 0\leq \theta\leq\frac{\pi}{2}.
\label{eq:kat}
\end{equation}
Moreover $\gamma^*_{\beta}$ is a strictly decreasing function of $\beta$ and in particular:
\begin{eqnarray}
&& \mbox{for $\pi\leq\beta\leq\beta_{cr}\;\;\;$ we have $\;\;\; 0.701\pi \, \approx \, \gamma^*_{\beta_{cr}}  \leq \gamma^*_{\beta} \leq \gamma^*_{\pi} \, \approx \, 0.867 \pi$} \nonumber \\[0.2cm]
&& \mbox{for $\beta_{cr}\leq\beta\leq 2\pi\;\;\;$ we have $\;\;\; 0.673 \pi \,\approx \, \gamma^*_{2\pi}  \leq \gamma^*_{\beta} \leq \gamma^*_{\beta_{cr}} \, \approx \, 0.701 \pi$.} \label{app}
\end{eqnarray}
\label{lemma8}
\end{lemma}
{\em Proof.} Inequality (\ref{eq:kat}) is written equivalently
\begin{equation}
\cot\frac{\gamma}{2} \geq \frac{\sin\theta}{\cos\theta +\frac{\alpha}{g(\beta,\theta)}},
\label{kat1}
\end{equation}
so what matters is the maximum of the function at the RHS of (\ref{kat1}). For each $0<\theta\leq \pi/2$ this function is strictly monotone as a function of $\beta$; this follows from Lemma \ref{lem:h} for $\beta_{cr}\leq\beta\leq 2\pi$
and from Lemma \ref{222aa} for $\pi\leq\beta\leq\beta_{cr}$.

The angle $\gamma^*_{\beta}\in (0,\pi)$ defined by
\[
\cot\frac{\gamma^*_{\beta}}{2} =  \max_{[0,\pi/2]} \frac{\sin\theta}{\cos\theta +\frac{\alpha}{g(\beta,\theta)}}
\]
is then a strictly increasing function of $\beta$. The approximate values in the statement have been obtained by numerical computations; see however Lemma \ref{lem:approx}. $\hfill\Box$

It would be nice to have good estimates on $\gamma^*_{\beta}$ without using a numerical solution of the differential equation (\ref{deg})
solved by $g(\theta)$.
This will be done for $\beta_{cr}\leq\beta\leq 2\pi$ by obtaining very good upper estimates on $g(\beta,\theta)$. We define
\[
\overline{g}(\beta , \theta) = a - \frac{a}{2(2a+1)}\theta^2 +\frac{a(4a^2+2a+3)}{24(2a+1)(4a^2+8a+3)}\theta^4 , \quad 0<\theta<\frac{\pi}{2},
\]
where $a$ is the largest solution of $a(1-a)=c_{\beta}$. We define the auxiliary quantity
$\gamma^{**}_{\beta}\in (0,\pi)$ by
\[
\cot\frac{\gamma^{**}_{\beta}}{2} =  \max_{[0,\pi/2]} \frac{\sin\theta}{\cos\theta +
\frac{\alpha}{\overline{g}(\beta,\theta)}}.
\]
\begin{lemma}
Let $\beta_{cr}\leq\beta\leq 2\pi$. Then  we have
\begin{eqnarray*}
\ia && g(\beta,\theta) \leq \overline{g}(\beta,\theta) \; , \qquad 0<\theta<\frac{\pi}{2},  \\[0.2cm]
\ib && \gamma^{**}_{\beta}\leq  \gamma^*_{\beta}.
\end{eqnarray*}
Actually we have (cf (\ref{app}))
\[
\gamma^{**}_{\beta_{cr}} \approx  0.700 \pi  \; ,  \qquad \gamma^{**}_{2\pi} \approx  0.672 \pi \, .
\]
\label{lem:approx}
\end{lemma}
{\em Proof.} We have $g(\beta,0)=\overline{g}(\beta,0)=\alpha$. Therefore, given that $g(\beta,\theta)$ satisfies
\begin{equation}
\label{g1}
\parder{g}{\theta} = -\frac{1}{\sin\theta}\Big( g^2 -g \cos\theta  + c_{\beta} \Big),
\end{equation}
it is enough to show that
\begin{equation}
\label{ug}
\parder{\overline g}{\theta} \geq -\frac{1}{\sin\theta}\Big( \overline{g}^2 -\overline{g} \cos\theta  + c_{\beta} \Big).
\end{equation}
The function $\overline{g}(\beta,\theta)$ is decreasing with respect to $\theta$, hence
\begin{eqnarray}
 \sin\theta \frac{d\overline{g}}{d\theta} +\overline{g}^2 &-& (\cos\theta)\overline{g} + c_{\beta} \nonumber \\ &\geq&
 \Big( \theta -\frac{\theta^3}{6}+ \frac{\theta^5}{120}  \Big)\frac{d\overline{g}}{d\theta}
 +\overline{g}^2 - \Big( 1-\frac{\theta^2}{2}+\frac{\theta^4}{24} \Big)\overline{g} + c_{\beta}\, .
\label{sto}
\end{eqnarray}
Now, a direct computation shows that the RHS of (\ref{sto}) is equal to
\[ \frac{a(1-a)\theta^6 [  16(2a+3)(2a+1)(22a^2+2a+3)-(12a^2+2a+3)(4a^2+2a+3)\theta^2 ]}
{2880(2a+1)^2(4a^2+8a+3)^2 } \]
\begin{eqnarray*}
&\geq& \frac{a(1-a)(12a^2+2a+3)(4a^2+2a+3)\theta^6 (16-\theta^2)}{2880(2a+1)^2(4a^2+8a+3)^2 } \\
&\geq &0.
\end{eqnarray*}
We note that in our argument we only used that $\alpha \in [1/2,1)$.

We now establish (i) for $\beta_{cr}<\beta\leq 2\pi$. The function 
\[
\overline{h}(\alpha,\theta)=\frac{\overline{g}(\beta,\theta)}{\alpha}
\]
(where, as usual, $\alpha$ is the largest solution of $\alpha(1-\alpha)=c_{\beta}<1/4$) is an upper solution to the initial value problem (\ref{ivp}). Hence applying (iii) of Lemma \ref{lem:h} we obtain the comparison.

To obtain (i) for $\beta=\beta_{cr}$ we use the monotonicity with respect to $\alpha$ of $h(\alpha,\theta)$. Passing to the limit $\alpha\to 1/2 +$ we conclude that
\[
H(\theta):= \lim_{\alpha \to 1/2 +}h(\alpha,\theta)\leq \overline{h}(\frac{1}{2},\theta)\leq 
2\overline{g}(\beta_{cr},\theta) \; , \qquad 0<\theta<\frac{\pi}{2}.
\]
The function $H(\theta)$ is then the maximal solution of the singular initial value problem (\ref{ivp}) and therefore coincides with the function $2g(\beta_{cr},\theta)$. This completes the proof of (i). Part (ii) then follows immediately from (i). $\hfill\Box$




\section{Proof of Theorem 1}
\label{sec:quads}

In this section we give the proofs of our theorems.
We start with a proposition that is fundamental in our argument and will be used repeatedly. We do not try to obtain the most general statement and for simplicity we restrict ourselves to assumptions that are sufficient for our purposes.

Let $U$ be a domain and assume that $\partial U=\Gamma\cup\Gamma_0$ where $\Gamma$ is Lipschitz continuous. We denote by $\vec{n}$ the exterior unit normal on $\Gamma$.
\begin{proposition}
Let $\phi\in H^1_{\rm loc}(U)$ be a positive function such that $\nabla\phi /\phi \in L^2(U)$ and $\nabla\phi /\phi$ has an $L^1$ trace on $\Gamma$ in the sense that $v \nabla\phi/\phi$ has an $L^1$ trace on $\partial U$ for every  $v\in C^{\infty}(\overline{U})$ that vanishes near $\Gamma_0$.
Then
\begin{equation}
\int_U |\nabla u|^2dx\, dy    \geq -\int_U\frac{\Delta\phi}{\phi}u^2dx\, dy    +\int_{\Gamma}  \frac{\nabla\phi}{\phi} \cdot\vec{n} u^2 dS
\label{lib}
\end{equation}
for all smooth functions $u$ which vanish near $\Gamma_0$ and $\Delta\phi$ is understood in the weak sense. \nl
If in particular there exists $c\in\R$ such that
\begin{equation}
-\Delta\phi\geq\frac{c}{d^2}\phi \; , 
\label{1}
\end{equation}
in the weak sense in $U$, where $d=\dist(x,\Gamma_0)$, then
\begin{equation}
\int_U |\nabla u|^2dx\,  dy  \geq c\int_{U}\frac{u^2}{d^2}dx\, dy   +\int_{\Gamma} u^2 \frac{\nabla\phi}{\phi} \cdot\vec{n} dS
\label{2}
\end{equation}
for all functions $u\in C^{\infty}(\overline{U})$ that vanish near $\Gamma_0$.
\label{lem:1}
\end{proposition}
{\em Proof.} Let $u$ be a function in $C^{\infty}(\overline{U})$ that vanishes near $\Gamma_0$. We denote $\vec{T}=-\nabla\phi/\phi$. Then
\begin{eqnarray*}
\int_U u^2 \diver \vec{T}\, dx\, dy &=& -2\int_U u\nabla u\cdot \vec{T}\, dx\, dy +\int_{\Gamma}u^2 \vec{T}\cdot \vec{n} \, dS  \\
&\leq &\int_U  |\vec{T}|^2u^2 dx\, dy +\int_U |\nabla u|^2 dx\, dy + \int_{\Gamma}u^2 \vec{T}\cdot \vec{n} \, dS\, ,
\end{eqnarray*}
that is
\[
\int_U |\nabla u|^2 dx\, dy \geq \int_U (\diver \vec{T} -|\vec{T}|^2) u^2 dx\, dy -\int_{\Gamma}  \vec{T}\cdot \vec{n}u^2 dS\, .
\]
Using assumption (\ref{1}) we obtain (\ref{2}). $\hfill\Box$






For $\beta\in (\pi,2\pi]$ we denote by $\Pi_{\beta}$ the class of all planar polygons which have precisely one non-convex vertex and the angle at that vertex is $\beta$.
Given a polygon in $\Pi_{\beta}$ we denote by $\gamma_+$ and $\gamma_-$ the angles at the vertices next to the non-convex vertex.

\begin{theorem}
Let $\beta\in (\pi,2\pi]$. Let $\Omega$ be a polygon in $\Pi_{\beta}$ with
\begin{equation}
\gamma_+,\gamma_- \leq \min \{ \gamma^*_{\beta} , \frac{3\pi -\beta}{2} \}
\label{cond:g1}
\end{equation}
where $\gamma^*_{\beta}\in (0,\pi)$ is defined by
\[
\cot\frac{\gamma^*_{\beta}}{2} =  \max_{[0,\pi/2]} \frac{\sin\theta}{\cos\theta +\frac{\alpha}{g(\beta,\theta)}}.
\]
Then the Hardy constant of $\Omega$ is $c_{\beta}$.
\label{new1}
\end{theorem}
{\em Proof.} We denote by $A_-$, $A_+$ the vertices next to the non-convex vertex $O$, so that $A_-$, $O$ and $A_+$ are consecutive vertices with respective angles
$\gamma_-$, $\beta$ and $\gamma_+$. We may assume that $O$ is the origin and that $A_+$ lies on the positive $x$-semiaxis.
We write the boundary $\partial\Omega$ as
\[
\partial\Omega = S_1 \cup S_2
\]
where $S_1 =OA_+ \cup OA_-$ and $S_2=\partial \Omega\setminus S_1$. We then define the equidistance curve
\[
\Gamma=\{ x\in\partial\Omega : \dist(x,S_1)= \dist(x,S_2) \}.
\]
Hence $\Gamma$ divides $\Omega$ into two sets $\Omega_1$ and $\Omega_2$, whose nearest boundary points belong in $S_1$ and $S_2$ respectively.
It is clear that $\Gamma$ can be parametrized by the polar angle $\theta\in [0,\beta]$.

The curve $\Gamma$ consists of line segments and parabola segments. Starting from $\theta=0$ we have line segments $L_1,\ldots,L_k$; then from
$\theta=\pi/2$ to $\theta=\beta -\pi/2$ we have parabola segments $P_1,\ldots,P_m$; and from $\theta=\beta -\pi/2$ to $\theta=\beta$ we have again line segments $L_1',\ldots,L_n'$.

Let $u\in C^{\infty}_c(\Omega)$ be given. Let $\vec{n}$ denote the unit normal along $\Gamma$ which is outward with respect to $\Omega_1$. Applying Proposition \ref{lem:1} with $\phi(x,y)=\psi_{\beta}(\theta)$, where $\theta$
is the polar angle of the point $(x,y)$, we obtain
\begin{equation}
\int_{\Omega_1}|\nabla u|^2 dx\, dy \geq c_{\beta}\int_{\Omega_1}\frac{u^2}{d^2}dx\, dy + 
\int_{\Gamma}\frac{\nabla\phi}{\phi}\cdot\vec{n} \, u^2 dS .
\label{pa}
\end{equation}
We next apply Proposition \ref{lem:1} on $\Omega_2$ for the function
$\phi_1(x,y)=d(x,y)^{\alpha}$ (we recall that $\alpha$ is the largest solution of $\alpha(1-\alpha)=c_{\beta}$). In $\Omega_2$ the function $d(x,y)$ coincides with the distance from $S_2$ and this implies that
\[
-\Delta d^{\alpha} \geq \alpha(1-\alpha)\frac{d^{\alpha}}{d^2}\; , \qquad \mbox{ on }\Omega_+\, .
\]
Applying Proposition \ref{lem:1} we obtain that 
\begin{eqnarray}
\int_{\Omega_+}|\nabla u|^2 dx\, dy &\geq& c\int_{\Omega_+}\frac{u^2}{d^2}dx\, dy  -\int_{\Gamma}\frac{\alpha\nabla d}{d}\cdot\vec{n}\, u^2 dS\, . 
\label{pa1}
\end{eqnarray}
Adding (\ref{pa}) and (\ref{pa1}) we conclude that 
\begin{equation}
\int_{\Omega}|\nabla u|^2 dx\, dy \geq c\int_{\Omega}\frac{u^2}{d^2}dx\, dy + 
\int_{\Gamma}\Big(\frac{\nabla\phi}{\phi}-\alpha\frac{\nabla d}{d}\Big)\cdot\vec{n}\, u^2 dS \, .
\label{nd}
\end{equation}
We emphasize that in the last integral the values of $\nabla\phi/\phi$ are obtained as limits from $\Omega_1$ and, more
importanmtly, those of $\nabla d/d$ are obtained as limits from $\Omega_2$.

It remains to prove that the line integral in (\ref{nd}) is non-negative. For this we shall consider the different segments of $\Gamma$. Due to the symmetry of our assumptions with respect to $\theta=\beta/2$ it is enough to establish the result for $0\leq \theta\leq \beta/2$.

(i) Let $L$ be one of the line segments $L_1,\ldots,L_k$. The points on this segment $L$ are equidistant from the side $OA_+$ and some side $E$ of $\partial\Omega \setminus (OA_+\cup OA_-)$.
Let $\gamma$ be the angle formed by the line $E$ and the $x$-axis so that the outward normal vector along $E$ is $(\sin\gamma,\cos\gamma)$ and $E$ has equation $x\cos\gamma +y\sin\gamma +c=0$ for some $c\in\R$. Elementary geometric considerations then give
$\gamma \in (-\pi/2,\pi)$.
Now, simple computations give
\begin{equation}
\Big(\frac{\nabla\phi}{\phi}-\alpha\frac{\nabla d}{d}\Big)\cdot\vec{n} =\frac{1}{d}\Big(g(\theta)\cos(\theta +\frac{\gamma}{2}) +
\alpha \cos(\frac{\gamma}{2})\Big) \; , \qquad \mbox{ on $L$.}
\label{f14}
\end{equation}
It remains to show that the RHS of (\ref{f14}) is non-negative for $0\leq\theta\leq\pi/2$. In the case $0<\gamma<\pi$  this is equivalent to showing that
\begin{equation}
 \cot\frac{\gamma}{2} \geq \frac{\sin\theta}{\cos\theta +\frac{\alpha}{g(\theta)}} , \quad 0\leq\theta\leq\frac{\pi}{2}.
\label{ch}
\end{equation}
This is true since $\gamma\leq\gamma_+\leq \gamma^*_{\beta}$.

In the case $-\pi/2 < \gamma \leq 0$ we have $\cos(\theta +\frac{\gamma}{2})\geq 0$ for all $0\leq\theta\leq\pi/2$ and the RHS is clearly non-negative.

(ii) Let $P$ be one of the parabola segments $P_1,\ldots,P_m$. The points on $P$ are equidistant from the origin $O$ and some side $E$ of $\partial\Omega \setminus (OA_+\cup OA_-)$. As in (i) above, 
let $\gamma$ be the angle formed by the line $E$ and the $x$-axis 
so that the outward normal vector along $E$ is $(\sin\gamma,\cos\gamma)$
and $E$ has equation $x\cos\gamma +y\sin\gamma +c=0$ for some $c\in\R$. Then $\gamma\in [\pi - \frac{\beta}{2},\pi ]$.
We note that the axis of the parabola has an asymptote at angle $\theta=\frac{3\pi}{2}-\gamma$. Indeed we shall prove the required
inequality for all $\theta\in [\frac{\pi}{2} , \frac{3\pi}{2}-\gamma] \supset [\frac{\pi}{2} , \frac{\beta}{2}]$.

Simple computations on $ P $ give
\begin{equation}
\Big(\frac{\nabla\phi}{\phi}-
\alpha\frac{\nabla d}{d}\Big)\cdot\vec{n} =\frac{1}{r\sqrt{2+2\sin(\theta+\gamma)}}\Big(f(\theta)\cos(\theta+\gamma) +\alpha
[1+\sin(\theta +\gamma )]\Big)   \; .
\label{fl6}
\end{equation}
Hence, noting that $\gamma\leq \gamma_+$, the result follows from Lemma \ref{lem:quad10}. This completes the proof.
$\hfill\Box$

{\bf\em Proof of Theorem 1.} This follows easily by approximating the convex set $K$ by a sequence of convex polygons and using Theorem \ref{new1}; see Fig.~\ref{fig:thm1}. $\hfill\Box$

{\bf Remark.} In case $\beta\leq\beta_{cr}$ we have $\gamma_{\beta}^*\leq \gamma_{\beta_{cr}}^* \approx 0.701\pi$ and therefore the condition $\gamma_+,\gamma_- \leq \min \{ \gamma^*_{\beta} , \frac{3\pi -\beta}{2} \}$
of Theorems \ref{thm1} and \ref{new1} takes the simpler form
\[
\gamma_+ ,  \gamma_- \leq  \gamma^*_{\beta}.
\]

If the convex set $K$ is unbounded and $\partial K$ does not intersect the boundary of $\Lambda_{\beta}$ then there is no need for any restriction. In particular
\begin{theorem}
Let $\Omega =K\cap\Lambda_{\beta}$ $K$ is an unbounded convex set and $\Lambda_{\beta}$ is a sector of angle $\beta\in (\pi,2\pi]$ whose vertex is inside $K$. Assume that the boundaries of $K$ and $\Lambda_{\beta}$ do not intersect. Then the Hardy constant of $\Omega$ is $c_{\beta}$, where $c_{\beta}$ is given by (\ref{tans1}), (\ref{tans2}).
\label{unbdd}
\end{theorem}
{\em Proof.} Let $u\in C^{\infty}_c(\Omega)$ be fixed. There exists a bounded convex set $K_1$ such that $\Omega_1:=K_1\cap S_{\beta}$ satisfies all the assumptions of Theorem \ref{thm1} and in addition
\[
{\rm dist}(x,\partial\Omega) =  {\rm dist}(x,\partial \Omega_1) \; , \qquad  x\in {\rm supp}(u)\, ;
\]
of course, $K_1$ depends on $u$. Applying Theorem \ref{thm1} to $\Omega_1$ we obtain the required Hardy inequality.
$\hfill\Box$

{\bf Remark.} Of course, one could state an intermediate result where the intersection $\partial K \cap \partial \Lambda_{\beta}$ is exactly one point forming an angle $\gamma$; in this the assumption $\gamma \leq \min \{ \gamma^*_{\beta} , \frac{3\pi -\beta}{2} \}$ should hold.











\section{Domains $E_{\beta,\gamma}$ with two non-convex angles}
\label{sec:ebg}

We reacall from the Introduction that
given angles $\beta$ and $\gamma$, we denote by $E_{\beta,\gamma}$ the domain shown in Fig.~\ref{fig:ebg1} in case
$\gamma<\pi$ and in Fig.~\ref{fig:ebg2} in case $\gamma >\pi$.
Its boundary $\partial E_{\beta,\gamma}$ consists of three parts $L_1$, $L_2$ and $L_3$. $L_2$ is a line segment and meets the halflines $L_3$ and $L_1$ at the origin O and the point $P(1,0)$ respectively.
We assume that $\beta+\gamma\leq 3\pi$ so that the halflines $L_1$ and $L_3$ do not intersect.
Without loss of generality we assume that $\beta\geq\gamma$ and since we are interested in the non-convex case, we assume that $\beta>\pi$.

{\bf\em Proof of Theorem \ref{thm2} part (i).} We denote by $\Gamma$ the curve
\[
\Gamma=\{ (x,y)\in E_{\beta,\gamma} : \dist((x,y), L_1) =\dist((x,y), L_2\cup L_3)\}.
\]
The curve $\Gamma$ divides $E_{\beta,\gamma}$ in two sets $E_-=\{ (x,y)\in E_{\beta,\gamma} : d(x,y)=\dist((x,y), L_2\cup L_3)\}$
and $E_+=\{ (x,y)\in E_{\beta,\gamma} : d(x,y)=\dist((x,y), L_1)\}$. We denote by $\vec{n}$ the unit normal along $\Gamma$ which is outward with respect to $E_-$.

Once again we shall use Proposition \ref{lem:1}.
We distinguish two cases: Case A, where $0\leq\gamma\leq \pi/2$ and Case B, where $\pi/2\leq\gamma\leq \pi$.

{\em Case A} ($0\leq\gamma\leq \pi/2$) We distinguish two subcases.

{\em Subcase Aa.}  $\beta+\gamma < 2\pi$. In this case $\Gamma$ consists of three parts: a line segment $\Gamma_1$ which bisects the angle at $P$; a parabola segment $\Gamma_2$, whose points are equidistant from the origin and the line $L_1$; and a halfline $\Gamma_3$ whose points are equidistant from $L_1$ and $L_3$. We parametrize $\Gamma$ by the polar angle $\theta$, so that $\Gamma_1=\{ 0\leq\theta\leq\frac{\pi}{2}\}$, $\Gamma_2=\{ \frac{\pi}{2}\leq\theta\leq\beta-\frac{\pi}{2}\}$,
and $\Gamma_3=\{  \beta-\frac{\pi}{2} \leq\theta < \frac{\beta+\pi-\gamma}{2}\}$.

Let $u\in C^{\infty}_c(E_{\beta,\gamma})$. We apply Proposition \ref{lem:1} with $U=E_-$, $\Gamma_0=L_2\cup L_3$ and for the function $\phi(x,y)=\psi(\theta)$, where $\psi=\psi_{\beta}$ and $\theta$
is the polar angle of $(x,y)$. We obtain that
\begin{equation}
\int_{E_-}|\nabla u|^2 dx\, dy \geq c_{\beta}\int_{E_-}\frac{u^2}{d^2}dx\, dy + 
\int_{\Gamma}\frac{\nabla\phi}{\phi}\cdot\vec{n}u^2 \,dS \, .
\label{pae}
\end{equation}
We next apply Proposition \ref{lem:1} to the domain $E_+$ and the function $\phi_1(x,y)=d(x,y)^{\alpha}$. We obtain that
\begin{equation}
\int_{E_+}|\nabla u|^2 dx\, dy \geq c_{\beta}\int_{E_+}\frac{u^2}{d^2}dx\, dy  -\alpha\int_{\Gamma}\frac{\nabla d}{d}\cdot\vec{n}u^2 dS\, .
\label{pa1e}
\end{equation}
Adding (\ref{pae}) and (\ref{pa1e}) we conclude that 
\begin{equation}
\int_{E_{\beta,\gamma}}|\nabla u|^2 dx\, dy \geq c_{\beta}\int_{E_{\beta,\gamma}}\frac{u^2}{d^2}dx\, dy + 
\int_{\Gamma}\Big(\frac{\nabla\phi}{\phi}-\alpha\frac{\nabla d}{d}\Big)\cdot\vec{n} \, u^2 dS \, .
\label{nde}
\end{equation}
We note that in the last integral the values of $\nabla\phi/\phi$ are obtained as limits from $E_-$ while those of $\nabla d/d$ are obtained as limits from $E_+$. 
It remains to prove that the last integral in (\ref{nde}) is non-negative. For this we shall consider the different parts of $\Gamma$.

(i) The segment $\Gamma_1$ ($0\leq\theta\leq\pi/2$). Simple computations give that 
\[
\frac{\nabla\phi}{\phi}-\alpha\frac{\nabla d}{d} =\frac{1}{d}\Big( g(\theta)\cos(\theta +\frac{\gamma}{2}) +
\alpha \cos(\frac{\gamma}{2})\Big)  \; , \qquad   0<\theta\leq\frac{\pi}{2}  \, ;
\]
this is non-negative by Lemma \ref{lemma8}, since $\gamma_{\beta}^*>\pi/2$.

(ii) The segment $\Gamma_2$ ($\pi/2\leq\theta\leq\beta-\pi/2$). In this case we have
\[
\Big(\frac{\nabla\phi}{\phi}-
\alpha\frac{\nabla d}{d}\Big)\cdot\vec{n} =\frac{1}{r\sqrt{2+2\sin(\theta+\gamma)}}\Big(f(\theta)\cos(\theta+\gamma) +\alpha
[1+\sin(\theta +\gamma )]\Big)  \;\; , 
\]
this is non-negative by Lemma \ref{lem:quad10}, since $\beta -\frac{\pi}{2} <\frac{3\pi}{2} -\gamma$.

(iii) The segment $\Gamma_3$ ($\beta-\frac{\pi}{2} \leq\theta <\frac{\beta+\pi-\gamma}{2} $). The line containing $\Gamma_3$ has equation
\[
x\cos(\frac{\beta-\gamma}{2}) + y\sin(\frac{\beta-\gamma}{2})  =\frac{\sin\gamma}{2\sin(\frac{\beta+\gamma}{2})},
\]
hence the outer (with respect to $E_-$) unit normal along $\Gamma_3$ is $(\cos(\frac{\beta-\gamma}{2}), \sin(\frac{\beta-\gamma}{2}))$.
Using the fact that $d=r\sin(\beta-\theta)$ on $\Gamma_3$, we have along $\Gamma_3$,
\begin{eqnarray*}
\Big(\frac{\nabla\phi}{\phi}-\alpha\frac{\nabla d}{d}\Big)\cdot\vec{n}&=&
[\frac{1}{r}\frac{\psi'(\theta)}{\psi(\theta)}(-\sin\theta , \cos\theta)  +\alpha\frac{(\sin\gamma ,\cos\gamma) }{d}] \\
&\cdot& (\cos(\frac{\beta-\gamma}{2}), \sin(\frac{\beta-\gamma}{2}))\\
&=&\frac{1}{r}\Big[ \frac{\psi'(\theta)}{\psi(\theta)}\sin(\frac{\beta-\gamma}{2}-\theta) +
\alpha \frac{\sin(\frac{\beta+\gamma}{2})}{\sin(\beta-\theta)} \Big]\\
&\geq& 0,
\end{eqnarray*}
since both terms in the last sum are non-negative (the first one, as the product of two non-positive terms).

{\em Subcase Ab.}  $\beta+\gamma \geq 2\pi$. In this case $\Gamma$ consists of only two parts $\Gamma_1$ and $\Gamma_2$, described exactly as in subcase Aa, the only difference being that the range of $\theta$ in $\Gamma_2$ is
$\frac{\pi}{2}\leq\theta<\frac{3\pi}{2} -\gamma$. This means that the parabola segment goes all the way to infinity. As before we have
\[
\Big(\frac{\nabla\phi}{\phi}-\alpha\frac{\nabla d}{d} \Big)\cdot \vec{n}=\frac{1}{r\sqrt{2+2\sin(\theta+\gamma)}}
\bigg(\frac{\psi'(\theta)}{\psi(\theta)}\cos(\theta+\gamma) +\alpha [1+\sin(\theta+\gamma)]\bigg)\\
\]
and the result follows again from Lemma \ref{lem:quad10}. This completes the proof in the case $0<\gamma\leq\pi/2$.

{\em Case B} ($\pi/2\leq\gamma\leq \pi$). On $E_-$ we again consider the function $\phi(x,y)=
\psi(\theta)$ and apply Lemma \ref{lem:1} as in the previous case. We fix a function $u\in C^{\infty}_c(E_{\beta,\gamma})$ and we obtain
\begin{equation}
\int_{E_-}|\nabla u|^2 dx\, dy \geq c_{\beta}\int_{E_-}\frac{u^2}{d^2}dx\, dy + 
\int_{\Gamma}(\frac{\nabla\phi}{\phi}\cdot\vec{n})u^2 dS \, .
\label{aaa1}
\end{equation}

In $E_+$ we consider a new orthonormal coordinate system with cartesian coordinates denoted by $(x_1,y_1)$ and polar coordinates denoted by $(r_1,\theta_1)$. The origin $O_1$ of this system is located on  the line $L_1$ and is such that the line $OO_1$ is perpendicular to $L_1$. The positive $x_1$ axis is then chosen so as to contain $L_1$ (diagram)
We note that this choice is such that
\begin{equation}
\label{113}
\mbox{the point on $\Gamma_1$ for which $\theta=\frac{\pi}{2}-\frac{\gamma}{2}$ satisfies also $\theta_1=\frac{\pi}{2}-\frac{\gamma}{2}$.}
\end{equation}

We apply Proposition \ref{lem:1} on $E_+$ with the function $\phi_1(x,y)=\psi(\theta_1)$. This function clearly satisfies
$-\Delta \phi_1 = c \, d^{-2}\phi_1$, hence we obtain
\begin{equation}
\int_{E+}|\nabla u|^2 dx\, dy \geq c \int_{E_+}\frac{u^2}{d^2}dx\, dy 
 -\int_{\Gamma}(\frac{\nabla\phi_1}{\phi_1}\cdot\vec{n})u^2  \,dS \; ,
\label{114}
\end{equation}
where, as before, $\vec{n}$ is the interior to $E_+$ unit normal along $\Gamma$.

Adding (\ref{aaa1}) and (\ref{114}) we conclude that
\begin{equation}
\int_{E_{\beta,\gamma}}|\nabla u|^2 dx\, dy \geq c_{\beta}\int_{E_{\beta,\gamma}}\frac{u^2}{d^2}dx\, dy + 
\int_{\Gamma}\Big(\frac{\nabla\phi}{\phi}-\frac{\nabla\phi_1}{\phi_1}\Big)\cdot\vec{n}\, u^2 dS \, .
\label{aaa3}
\end{equation}
The rest of the proof is devoted to showing that the last integral
in (\ref{aaa3}) is non-negative. 

As in the case $0<\gamma\leq \pi/2$, we need to distinguish two subcases: Subcase Ba, where $\beta+\gamma < 2\pi$, and Subcase Bb, where $\beta+\gamma\geq 2\pi$.

{\em Subcase Ba.} $\beta+\gamma < 2\pi$.
The curve $\Gamma$ consists of three parts: a line segment $\Gamma_1$ which bisects the angle at $P$; a (part of a) parabola $\Gamma_2$, whose points are equidistant from the origin and the line $L_1$; and a halfline $\Gamma_3$ whose points are equidistant from $L_1$ and $L_3$. As before, we consider separetely each segment and we parametrize $\Gamma$ by the polar angle $\theta$ so that
\[
\Gamma_1=\{ \theta\in\Gamma  : \,
0\leq\theta\leq\frac{\pi}{2}\} \; , \;\; \Gamma_2=\{ \frac{\pi}{2}\leq\theta\leq\beta-\frac{\pi}{2}\}
\; , \]
\[ \Gamma_3=\{\beta-\frac{\pi}{2} \leq\theta < \frac{\beta+\pi-\gamma}{2}\}. \]
(i) The segment $\Gamma_1$ ($0\leq\theta\leq\pi/2$). We have
\[
\frac{\nabla\phi}{\phi}\cdot\vec{n} = \frac{\psi'(\theta)}{r\psi(\theta)}\cos(\theta+\frac{\gamma}{2}) \; , \quad
\mbox{ on } \Gamma_1.
\]
and similarly
\[
\frac{\nabla\phi_1}{\phi_1}\cdot\vec{n} = -\frac{\psi'(\theta_1)}{r_1 \psi(\theta_1)}\cos(\theta_1-\frac{\gamma}{2}) \; , \quad
\mbox{ on }  \Gamma_1.
\]
Since $r_1\sin\theta_1 =r\sin\theta$ along $\Gamma_1$, it is enough to prove the inequality
\begin{equation}
g(\theta)\cos(\theta+\frac{\gamma}{2}) +g(\theta_1)\cos(\theta_1-\frac{\gamma}{2})\geq 0 \; , \qquad 
0\leq \theta\leq \frac{\pi}{2}.
\label{222}
\end{equation}
This has been proved in \cite{BT}; we include a proof here for the sake of completeness.
Recalling (\ref{113}) and  applying the sine law we obtain that along $\Gamma_1$ the polar angles $\theta$ and $\theta_1$ are related by
\begin{equation}
\cot \theta_1 =-\cos\gamma \cot\theta +\sin\gamma \; .
\label{223}
\end{equation}
{\bf Claim.} There holds
\begin{equation}
\label{claim1}
 \theta_1\geq  \theta +\gamma -\pi \; , \qquad \mbox{ on }\Gamma_1 \, .
\end{equation}
{\em Proof of Claim.} We fix $\theta\in [0,\pi/2]$ and 
the corresponding $\theta_1=\theta_1(\theta)$. If $\theta+\gamma -\pi\leq 0$, then (\ref{claim1}) is obviously true, so we assume that $\theta+\gamma -\pi\geq 0$. Since $0\leq\theta+\gamma -\pi\leq \pi/2$ and $0\leq\theta_1\leq\pi/2$, (\ref{claim1}) is written equivalently $\cot\theta_1\leq\cot(\theta +\gamma -\pi)$; thus, recalling (\ref{223}), we conclude that to prove the claim it is enough to show that
\[
 -\cos\gamma\cot\theta +\sin\gamma\leq\cot(\theta+\gamma) \; , \quad
\pi-\gamma\leq\theta\leq\frac{\pi}{2},
\]
or, equivalently (since $\pi\leq\theta+\gamma\leq 3\pi/2$),
\begin{equation}
-\cos\gamma\cot^2\theta +(-\cos\gamma\cot\gamma -\cot\gamma+\sin\gamma)\cot\theta + 1+\cos\gamma \geq 0\; , 
\pi-\gamma\leq\theta\leq\frac{\pi}{2}.
\label{eq:new}
\end{equation}
The left-hand side of (\ref{eq:new}) is an increasing function of $\cot\theta$ and therefore takes its least value at $\cot\theta =0$.
Hence the claim is proved.

For $0\leq\theta\leq\pi/2-\gamma/2$ (\ref{222}) is true since all terms in the left-hand side are non-negative. So let
$\pi/2-\gamma/2\leq\theta\leq\pi/2$ and $\theta_1=\theta_1(\theta)$. From (\ref{223}) we find that
\begin{eqnarray*}
\frac{d\theta_1}{d\theta}-1&=&-\frac{\cos\gamma (1+\cot^2\theta)+1+\cot^2\theta_1}{1+\cot^2\theta_1} \\
&=& -\frac{ 1+\sin^2\gamma +\cos\gamma -2\sin\gamma\cos\gamma\cot\theta +\cos\gamma(1+\cos\gamma)\cot^2\theta}{1+\cot^2\theta_1}.
\end{eqnarray*}
The function
\[
h(x):= 1+\sin^2\gamma +\cos\gamma -2\sin\gamma\cos\gamma x +\cos\gamma(1+\cos\gamma)x^2
\]
is a concave function of $x$. We will establish the positivity of $h(\cot\theta)$ for $\pi/2-\gamma/2
\leq \theta\leq\pi/2$. For this it is enough to establish the positivity at the endpoints. At $\theta=\pi/2$ positivity is obvious, whereas
\[
h(\tan(\frac{\gamma}{2}))=1+\sin^2\gamma+\cos\gamma -2\cos\gamma\sin^2\frac{\gamma}{2}\geq 0.
\]
From (\ref{113}) we conclude that $\theta_1\leq\theta$ for $\pi/2-\gamma/2\leq\theta\leq\pi/2$. Now, it was proved in \cite[Lemma 4]{BT} that the function $g$ is decreasing. Hence for $\pi/2-\gamma/2\leq\theta\leq\pi/2$ we have,
\begin{eqnarray*}
g(\theta)\cos(\theta+\frac{\gamma}{2}) +g(\theta_1)\cos(\theta_1-\frac{\gamma}{2})
&\geq& g(\theta) [\cos(\theta+\frac{\gamma}{2}) +\cos(\theta_1-\frac{\gamma}{2})]\\
&=&2g(\theta)\cos(\frac{\theta+\theta_1}{2}) \cos (\frac{\theta-\theta_1+\gamma}{2})\\
&\geq& 0,
\end{eqnarray*}
where for the last inequality we made use of the claim. Hence (\ref{222}) has been proved.

(ii) The segment $\Gamma_2$ ($\frac{\pi}{2}\leq\theta\leq\beta-\frac{\pi}{2}$). After some computations we obtain that
\begin{eqnarray*}
\Big(\frac{\nabla\phi}{\phi}-\frac{\nabla\phi_1}{\phi_1}\Big)\cdot\vec{n} =
\frac{1}{r\sqrt{2+2\sin(\theta+\gamma)}}
&\bigg\{&f(\theta)\cos(\theta+\gamma) \\
&-& f(\theta_1)\sin\theta_1[\sin(\theta_1-\theta -\gamma)-\cos\theta_1]\bigg\} ,
\end{eqnarray*}
where $\theta$ and $\theta_1$ are related by $\cot\theta_1 =-\cos(\theta+\gamma)$. The result then follows
by applying \cite[Lemma 6]{BT}.

(iii) The segment $\Gamma_3$ ($\beta-\frac{\pi}{2}\leq\theta< \frac{\beta+\pi-\gamma}{2}$). Simple computations yield that along $\Gamma_3$ we have
\begin{equation}
\label{fs}
\Big(\frac{\nabla\phi}{\phi}-\frac{\nabla\phi_1}{\phi_1}\Big)\cdot\vec{n} =\frac{\psi'(\theta)}{r\psi(\theta)}\sin(\frac{\beta-\gamma}{2}-\theta)
+\frac{\psi'(\theta_1)}{r_1\psi(\theta_1)}\sin(\frac{\beta+\gamma}{2}-\theta_1).
\end{equation}
The first summand in the right-hand side of (\ref{fs}) is non-negative since $\psi'(\theta)$ and $\sin(\frac{\beta-\gamma}{2}-\theta)$ are non-positive in the given range of $\theta$. Moreover, two applications of the sine law yield that along $\Gamma_3$ the coordinates $(r,\theta)$ and $(r_1,\theta_1)$ are related by
\[
r_1\sin\theta_1=r\sin(\beta-\theta) \;\; , \qquad \tan\theta_1 =-\frac{\sin(\beta-\theta)}{\cos(\theta+\gamma)}.
\]
It follows in particular that $0\leq\theta_1\leq \pi/2$, and hence $\pi/4 \leq\frac{\beta+\gamma}{2}-\theta_1 \leq \pi$. Hence the second summand in the right-hand side of (\ref{fs}) is also non-negative, completing the proof in this case.

{\em Subcase B2.} $\beta+\gamma\geq 2\pi$. In this case $\Gamma$ consists only of two parts $\Gamma_1$ and $\Gamma_2$, described as in Case B1. The only difference is that the range of $\theta$ in $\Gamma_2$ now is $\frac{\pi}{2}\leq\theta<\frac{3\pi}{2} -\gamma$; the result follows as before. This completes the proof of the theorem.
$\hfill\Box$

{\bf\em Proof of Theorem \ref{thm2} part (ii).} We set for simplicity $\psi=\psi_{\beta+\gamma-\pi}$. 
We divide $E_{\beta,\gamma}$ in three parts $E_1$, $E_2$ and $E_3$ as in the diagram, and denote $L_i=(\partial E_i)\cap\partial E_{\beta,\gamma}$. We also set $\Gamma_i=\{(i,y) : y\geq 0\}$, $i=0,1$, the halflines that are the common boundaries of the $E_j$'s.  We first apply Proposition \ref{lem:1} to the domain $E_1$.
For this we introduce polar coordinates $(r_1,\theta_1)$ centered at $P$, so that the positive $x_1$ axis coincides with the halfline $L_1$. Let $u\in C^{\infty}_c(E_{\beta,\gamma})$ be fixed.
Applying Proposition \ref{lem:1} with $\phi(x,y)=\psi(\theta_1)$ we obtain
\begin{equation}
\int_{E_1}|\nabla u|^2dx\, dy \geq c_{\beta+\gamma-\pi}\int_{E_1}\frac{u^2}{d^2}dx\, dy + \frac{\psi'(\gamma-\frac{\pi}{2})}{\psi(\gamma-\frac{\pi}{2})}\int_{\Gamma_1}\frac{u^2}{y}dy \; .
\label{11a}
\end{equation}
On $E_3$ we use the standard polar coordiantes $(r,\theta)$ and the function $\phi(x,y)=\psi(\beta-\theta)$. We obtain
\begin{equation}
\int_{E_3}|\nabla u|^2dx\, dy \geq c_{\beta+\gamma-\pi}\int_{E_3}\frac{u^2}{d^2}dx\, dy + \frac{\psi'(\beta-\frac{\pi}{2})}{\psi(\beta-\frac{\pi}{2})}\int_{\Gamma_0}\frac{u^2}{y}dy \; .
\label{11b}
\end{equation}
Without loss of generality we assume that $\beta\geq\gamma$ and we therefore have
\[
\frac{\psi'(\gamma-\frac{\pi}{2})}{\psi(\gamma-\frac{\pi}{2})}=-\frac{\psi'(\beta-\frac{\pi}{2})}{\psi(\beta-\frac{\pi}{2})}\geq 0. 
\]
Now, we have $u(1,y)^2 -u(0,y)^2 =2\int_0^1 uu_xdx$, hence, using also the 1-dimensional Hardy inequality we have for any $\epsilon>0$,
\begin{eqnarray*}
\int_{\Gamma_0}\frac{u^2}{y}dy - \int_{\Gamma_1}\frac{u^2}{y}dy&\leq& \epsilon\int_{E_2}\frac{u^2}{y^2}dx\, dy +\frac{1}{\epsilon}\int_{E_2}u_x^2dx\, dy\\
&\leq& (\epsilon-\frac{1}{4\epsilon})\int_{E_2}\frac{u^2}{y^2}dx\, dy  + \frac{1}{\epsilon}\int_{E_2} u_y^2 dx\, dy +\frac{1}{\epsilon}\int_{E_2}u_x^2dx\, dy
\end{eqnarray*}
and therefore
\begin{equation}
\int_{E_2}|\nabla u|^2 dx\, dy \geq \Big(\frac{1}{4}-\epsilon^2\Big) \int_{E_2}\frac{u^2}{y^2}dx\, dy +\epsilon
\int_{\Gamma_0}\frac{u^2}{y}dy - \epsilon\int_{\Gamma_1}\frac{u^2}{y}dy\, .
\label{11c}
\end{equation}
This is also true for $\epsilon=0$. We choose $\epsilon=\psi'(\gamma-\frac{\pi}{2})/\psi(\gamma-\frac{\pi}{2})$ 
and we note that by (\ref{oa}) we have
\[
c_{\beta+\gamma-\pi}\leq \frac{1}{4}-
 c_{\beta+\gamma-\pi}\tan^2\big( \sqrt{c_{\beta+\gamma-\pi}}\frac{\beta-\gamma}{2}\big)
 = \frac{1}{4}- \Big(\frac{\psi'(\gamma-\frac{\pi}{2})}{\psi(\gamma-\frac{\pi}{2})}\Big)^2 
 = \frac{1}{4} -\epsilon^2 \, .
\]
Adding (\ref{11a}), (\ref{11b}) and (\ref{11c}) we obtain the inequalities in all cases.

We now prove the sharpness of the constant. Let $C$ denote the best Hardy constant for $E_{\beta,\gamma}$.
We extend the halflines $L_1$ and $L_3$ until they meet at a point $A$, and we call $D_0$ the resulting infinite sector, whose angle is $\beta+\gamma-\pi$. We introduce a family of domains $D_{\epsilon}$ that are obtained from $E_{\beta,\gamma}$ by moving $L_2$ parallel to itself towards $A$ so that it is a distance $\epsilon$ from $A$. All these domains $D_{\epsilon}$ have the same Hardy constant as $E_{\beta,\gamma}$.
Let $d_{\epsilon}(x)={\rm dist}(x,\partial D_{\epsilon})$ and $d_0(x)={\rm dist}(x,\partial D_0)$. Then clearly $d_{\epsilon}(x)\to d_0(x)$ for all $x\in D_0$.

Let $u\in C^{\infty}_c(D_0)$ vanish near $\Gamma_0$.
This can be used as a test function for the Hardy inequality in $D_{\epsilon}$, therefore we have
\[
\int_{D_{\epsilon}}|\nabla u|^2 dx\, dy \geq C\int_{D_{\epsilon}} \frac{u^2}{d_{\epsilon}^2}dx\, dy \, ,
\]
which can be written equivalently
\[
\int_{D_0}|\nabla u|^2 dx\, dy \geq C\int_{D_0} \frac{u^2}{d_{\epsilon}^2}dx\, dy \, .
\]
Passing to the limit $\epsilon\to 0$ we therefore obtain
\[
\int_{D_0}|\nabla u|^2 dx\, dy \geq C\int_{D_0} \frac{u^2}{d_0^2}dx\, dy \, .
\]
Since the best Hardy constant of $D_0$ is $c_{\beta+\gamma-\pi}$, we conclude that $C\leq c_{\beta+\gamma-\pi}$, which establishes the sharpness. $\hfill\Box$





\section{A Dirichlet - Neumann Hardy inequality}

We finally prove Theorem \ref{thm3}.

{\bf\em Proof of Theorem \ref{thm3}.} Let $u\in C^{\infty}(\overline{D_{\beta}})$. Applying Proposition \ref{lem:1} for $\phi(x,y)=\psi(\theta)$ we have
\begin{eqnarray*}
\int_{D_{\beta}} |\nabla u|^2dx\, dy    &\geq&  -\int_{D_{\beta}}\frac{\Delta\phi}{\phi}u^2dx\, dy    +\int_{\Gamma}  \frac{\nabla\phi}{\phi} \cdot\vec{n} u^2 dS \\
&=& c_{\beta}\int_{D_{\beta}}\frac{u^2}{d^2}dx\, dy    +\int_{\Gamma}  \frac{\nabla\phi}{\phi} \cdot\vec{n} u^2 dS.
\end{eqnarray*}
A direct computation gives that along $\Gamma$ we have
\[
\frac{\nabla\phi}{\phi} \cdot\vec{n} = - \frac{r'(\theta)}{r(\theta) \sqrt{ r(\theta)^2 +r'(\theta)^2}}\cdot
\frac{\psi'(\theta)}{\psi(\theta)},
\]
which establishes the inequality. The fact that $c_{\beta}$ is sharp follows by comparing with the corresponding Dirichlet problem. $\hfill\Box$





\end{document}